\documentclass[fleqn,10pt]{article}
\usepackage{amsmath}
\usepackage{epsf}
\usepackage{epsfig}
\usepackage{subfigure}
\usepackage{a4}
\usepackage[margin=3.5cm]{geometry}
\usepackage{amssymb}
\usepackage{algorithm}
\usepackage{algorithmic}
\usepackage{epstopdf}
\def\kme{k_{\rm meso}}

\newtheorem{theorem}{Theorem}

\newtheorem{corollary}{Corollary}

\title{{
On the Reaction Diffusion Master Equation in the Microscopic Limit}}
\author{{\small 
      Stefan Hellander$^{\mbox{\tiny 1}}$, Andreas Hellander$^{\mbox{\tiny 2}}$, Linda Petzold$^{\mbox{\tiny 2}}$
       }
}
\date{}

\begin{document}

\maketitle

\vspace{-10pt}


\begin{center}
{\footnotesize\em 
$^{\mbox{\tiny\rm 1}}$Division of Scientific Computing,
Department of Information Technology \\
Uppsala University, P.~O.~Box 337, SE-75105 Uppsala, Sweden \\
email: {\tt stefan.hellander@it.uu.se} \\
$^{\mbox{\tiny\rm 2}}$Department of Computer Science,\\
University of California, Santa Barbara, CA 93106-5070 Santa Barbara, USA.\\
email: {\tt andreash@cs.ucsb.edu, petzold@cs.ucsb.edu}\\
[3pt]}
\end{center}

\begin{abstract}

Stochastic modeling of reaction-diffusion kinetics has emerged as a powerful theoretical tool in the study of biochemical reaction networks. Two frequently employed models are the particle-tracking Smoluchowski 
framework and the on-lattice Reaction-Diffusion Master Equation (RDME) framework. As the mesh size goes from coarse to fine, the RDME initially becomes more accurate. However, recent developments have shown that it will become increasingly inaccurate compared to the Smoluchowski model as the lattice spacing becomes very fine.
In this paper we give a new, general and simple argument for why the RDME breaks down. Our analysis reveals a hard limit on the voxel size for which no local RDME can agree with the Smoluchowski model.  

\end{abstract}
\section{Introduction}
A prevalent view in molecular systems biology is that the noise in cellular reaction networks, arising intrinsically from low copy numbers of macromolecules, can have a substantial impact on function \cite{Elowitz:2002fj,FaEl}. Two frequently used models for simulating stochastic reaction-diffusion systems are the reaction-diffusion master equation (RDME) \cite{VanKampen,Gardiner} and the Smoluchowski model \cite{Smol}, which we will refer to as the mesoscopic and microscopic models, respectively. In the RDME the computational domain is divided into voxels. Diffusion is simulated as discrete jumps between adjacent vertices in the mesh. Reactions can occur locally when reactants are present in the same voxel. The RDME is attractive from a computational perspective; it is the logical extension of spatially homogenous simulations based on the Gillespie algorithm \cite{Gillespie}, and keeps track of the location of molecules only up to the resolution of the mesh, hence allowing for coarse-graining.
Publicly available software packages for mesoscopic simulations based on the RDME include mesoRD \cite{mesoRD}, SmartCell \cite{SmartCell} and URDME \cite{urdme11}.   

On a finer modeling level, the Smoluchowski model treats diffusion and reaction in continuous space, with molecules explicitly represented as spheres with a certain interaction radius. As such, it is an example of a model commonly referred to as particle-tracking. Software based on the Smoluchowski model are Smoldyn \cite{And,AnAdBrAr10, AnB}, MCell \cite{MCell08} and Green's function reaction dynamics (GFRD) \cite{ZoWo5a,ZoWo5b}. Smoldyn and MCell both discretize time whereas GFRD is continuous in time, making it more accurate at fine scales. 

A well known property of the mesoscopic model is that it converges to the classical reaction-diffusion partial differential equation in the macroscopic limit. For a system approaching the microscopic regime, it is tempting to think of the RDME as a better and better approximation to the Smoluchowski model for finer and finer mesh resolutions. This picture is misleading, as in fact, it has been shown that as the size of the voxels in an infinite 3D domain decreases, all bimolecular reactions are eventually lost in the mesoscopic model \cite{Isaacson2}. Recent work using GFRD has demonstrated that fast, microscopic rebinding events can substantially affect the macroscopic properties of a biochemical signal cascade when reactions are highly diffusion limited \cite{TaTNWo10}. To accurately simulate such systems requires a fine spatial resolution. On these scales, the conventional RDME may be too inaccurate to capture even the qualitative behavior predicted by the microscopic model \cite{FBSE10}.         

The fact that the conventional mesoscopic model becomes inaccurate as we approach the microscopic level is not surprising, as we are moving out of the domain of validity for which it was derived. However it can pose a real practical problem, as it is hard to know \emph{a priori} if a simulation with the RDME will yield useful or misleading results. This is especially true for biochemical models with multiscale properties, which are frequently encountered in molecular biology. Simply resorting to e.g. GFRD whenever in doubt is currently not feasible in general due to the high computational cost for systems with many particles. A natural approach to remedy this problem is to try to extend the domain of validity of the RDME as the mesh size tends to zero. Isaacson \cite{Isaacson2} suggests that one way of doing this would be to let the association rate constants depend explicitly on the meshsize. Recently, approaches to make such corrections to the RDME have been proposed \cite{FBSE10,ErbanChapman}.   

In this paper we show, by a simple argument, that below a certain critical size of the mesh it will be impossible to make the RDME consistent with the Smoluchowski model by local modifications to the bimolecular rates. Our result is valid for finite domains in 2D and 3D, and for mesh-dependent mesoscopic association rates. For a simple model problem we derive the correct mesoscopic association rates in 2D and 3D and compute the critical size of the mesh, and show that it can be considerably larger than the interaction radii of the molecules. In this light we review and discuss recent work in which the RDME has been modified in different ways in order to better agree with the microscale model for very small voxel sizes.    

\section{Background}
\subsection{The Reaction-Diffusion Master Equation}
In the mesoscopic model, the computational domain is divided into non-overlapping voxels. The state of the system is the discrete number of molecules of each biochemical species in each voxel in the mesh. Inside voxels, the molecules react according to a pre-specified set of rules. In a small time interval $dt$, a bimolecular reaction $A+B \rightarrow C$, for example, occurs with probability $k_a a bdt/V$, where $k_a$ is the mesoscopic rate constant for the reaction, $a$ and $b$ are the copy numbers of $A$ and $B$ in that voxel, and $V$ is the volume of the voxel. Diffusion is modeled as jumps between adjacent voxels. For a Cartesian mesh with mesh spacing $h$, the rate for a diffusive jump  is given by $D/h^2$, where $D$ is the diffusion constant. The time evolution of the whole system is described as a continuous-time discrete-space Markov process and  the probability density function (PDF) of the system evolves according to the RDME \cite{VanKampen, Gardiner}. Solving for the PDF directly is impossible in general due to the extremely large statespace, but realizations of the process can be generated efficiently using kinetic Monte Carlo methodology \cite{ElE}.

\subsection{The Smoluchowski model}
In the Smoluchowski model two molecules $A$ and $B$ are assumed to move by Brownian motion with diffusion constants $D_A$ and $D_B$,  and react with a certain probability at a distance determined by the sum of their reaction radii, $\rho_A$ and $\rho_B$. For two molecules $A$ and $B$, the probability of a bimolecular reaction is governed by the Smoluchowski equation. Given that the molecules start at positions ${x_A}^0$ and ${x_B}^0$ at time $t_0$, the equation for the PDF $p$ of the new relative position $\mathbf{r}=x_A-x_B$ (in a spherical coordinate system $\mathbf{r} = (r,\theta,\phi)$), is given by $\partial_t \mathbf{r} = D\Delta \mathbf{r}$ with initial condition $p\left(\mathbf{r},t_0|\mathbf{r}_0,t_0\right) = \delta\left(\mathbf{r}-\mathbf{r}_0\right)$ and boundary conditions
\begin{align*}
&\lim_{|\mathbf{r}|\to\infty} p\left(\mathbf{r},t|\mathbf{r}_0,t_0\right) = 0,\\
&4\pi\rho D\left.\frac{\partial p}{\partial r}\right|_{r=\sigma} = k_r p\left(\mathbf{r},t|\mathbf{r}_0,t_0\right)|_{r=\sigma}
\end{align*}
where $\rho = \rho_A+\rho_B$, $D = D_A+D_B$ and $k_r$ is the microscopic association rate. It can be shown that a weighted mean of the positions given by $\mathbf{R} = \sqrt{D_B/D_A}x_A+\sqrt{D_A/D_B}x_B$ will be normally distributed \cite{ZoWo5a}, and by sampling a new $\mathbf{r}$ and $\mathbf{R}$ we obtain the new positions of the molecules at some time $t$.

An efficient method for simulating systems of molecules is the Green's function reaction dynamics (GFRD) \cite{ZoWo5a,ZoWo5b} method. In GFRD a system of molecules is decomposed into pairs of molecules by choosing a time step such that any given molecule is unlikely to react with more than one other molecule during that time step. This reduces the full problem to solving the Smoluchowski equation for pairs of molecules. 

\section{Results}
\subsection{Breakdown of the mesoscopic model}
Recent work has demonstrated that the RDME breaks down in the limit of infinitesimal voxels. Isaacson \cite{Isaacson2} shows that the probability for the occurrence of bimolecular reactions vanishes with decreasing voxel size for molecules on a lattice in an infinite 3D domain. The study is restricted to the case where the mesoscopic reaction rates are not dependent on the size of the voxels. Here we present a new and intuitive way to understand the degeneration of the mesoscopic model in finite domains in 2D and 3D, and with mesoscopic reaction rates that may or may not depend explicitly on the mesh. Our analysis will also provide additional insight into why and when this breakdown occurs. 

To see why the RDME model cannot work for very small voxels, it is illustrative to consider the simple process of bimolecular association 
\begin{align} 
A+B\xrightarrow{k_a} C.
\label{eq:dimer}
\end{align}
In 3D, the diffusion-controlled mesoscopic reaction rate $k_a$ is often taken to be
\begin{equation}
k_a = \frac{4\pi\rho D k_r}{4\pi\rho D + k_r}.
\label{eq:mesomicro}
\end{equation}
This expression will hereafter be referred to as the conventional mesoscopic rate constant. It is valid for large enough voxels, and in 2D no analogous expression is well-defined. 

Consider the reaction \eqref{eq:dimer} and, in a Cartesian coordinate system, let a molecule of $A$ be stationary in some voxel $V_A$, i.e. $D_A = 0$, and let $B$ diffuse with diffusion constant $D_B$ in a square or cube with side length $L$ and with periodic boundary conditions. Let $\tau_{\mathrm{meso}}$ be the average time until the two molecules react in the discretized, mesoscopic model, let $\tau_J$ denote the average time for a diffusive jump and let $\kme$ be the mesoscopic reaction rate. Once the $B$ molecule is in the voxel $V_A$, the two molecules will react with probability $\kme/(\kme+\tau_J)$, or the $B$ molecule diffuses away one voxel. Thus the two molecules will be on average $M = (\kme+\tau_J)/\kme$ times in the same voxel before they react. Let $\tau_D$ denote the average time until the two molecules are in the same voxel for the first time. The average time for an event (reaction or diffusion) given that the two molecules are in the same voxel is $\tau_e=1/(\kme+\tau_J)$, and thus the average time for the $B$ molecule to diffuse away one voxel and then back is $\tau_e+\tau_D^1$, where $\tau_D^1$ is the average time until the $B$ molecule reaches $V_A$ given that it starts one voxel away. This will be repeated on average $M-1$ times until the molecules react. We can now write $\tau_{\mathrm{meso}}$ as
\begin{align}
\tau_{\mathrm{meso}} &= \tau_D+(M-1)\left(\tau_e+\tau_D^1\right)+\tau_e \nonumber\\ &=\tau_D+\kme^{-1}(1+N^1_{\mathrm{steps}})=:\tau_D+\tau_r.
\label{eq:taumeso}
\end{align}
where $N^1_{\mathrm{steps}}=\tau_J\tau_D^1$. The following Theorem was proven in \cite{Montroll68}.
\begin{theorem}
Assume that the molecule $B$ has a uniformly distributed random starting position $x_B$ on the lattice, $x_B$ not equal to $x_A$, and that the molecules can move to nearest neighbors only (as in the RDME). Then as $N\to \infty$, the following holds:
\begin{align*}
&N_{\mathrm{steps}} \sim \pi^{-1}N\log(N)+0.1951N, \, N^1_{\mathrm{steps}} \sim N\, \mathrm{(2D)}\\
&N_{\mathrm{steps}} \sim 1.5164N, \, N^1_{\mathrm{steps}} \sim N, \, \mathrm{(3D)}.
\end{align*}  
where $N_{\mathrm{steps}}$ is the average number of steps until $x_B=x_A$ for the first time, $N_{\mathrm{steps}}^1=\tau_J\tau_D^1$ is the average number of steps until $x_B=x_A$ given that $A$ and $B$ start one voxel apart and $N$ is the number of voxels in the domain.
\label{theoremNsteps}
\end{theorem}
\begin{corollary}
Let $\tau_D$ be the time until $A$ and $B$ are in the same voxel for the first time. Then
\begin{align*}
&\tau_D \sim \frac{L^2}{2\pi D}\log\left(\frac{L}{h}\right)+0.1951\frac{L^2}{4D} \quad\text{as } h\to 0 & \quad \mathrm{(2D)},\\
&\tau_D \sim 1.5164\frac{L^3}{6Dh} \quad\text{as } h\to 0 & \quad \mathrm{(3D)}.
\end{align*}
where $h$ is the voxel size.
\label{cor1}
\end{corollary}
This follows immediately from Theorem \ref{theoremNsteps} and the fact that $\tau_D=\frac{N-1}{N}N_{\mathrm{steps}}\tau_J^{-1}$, where $\tau_J=2dD/h^2$ and $d$ is the dimension.

From Corollary \ref{cor1} and Equation \eqref{eq:taumeso} we conclude that for a sufficiently small voxel size in the discrete space model, we will have $ \tau_{\mathrm{meso}}>\tau_D>\tau_{\mathrm{micro}}$ (where $\tau_{\mathrm{micro}}$ is the average time until the molecules react in the microscopic model), \emph{for any choice} of the mesoscopic rate constant, since $\tau_r>0$ for all $\kme>0$. Eventually, as $h \to 0$, no bimolecular reactions will occur since molecules can only react when they are in the same voxel and $\tau_D\to\infty$, hence $\tau_{\mathrm{meso}} \to \infty$. 

Note that the reason for this effect is not that the diffusion process is inaccurately described at these length scales. The problem lies in the assumptions that molecules can react only after having diffused to the same voxel.
\subsection{No local correction to the association rates can make the RDME consistent with the microscopic model.}
\begin{figure}[h]
\centering
\subfigure[\, 2D]{\includegraphics[scale=0.5]{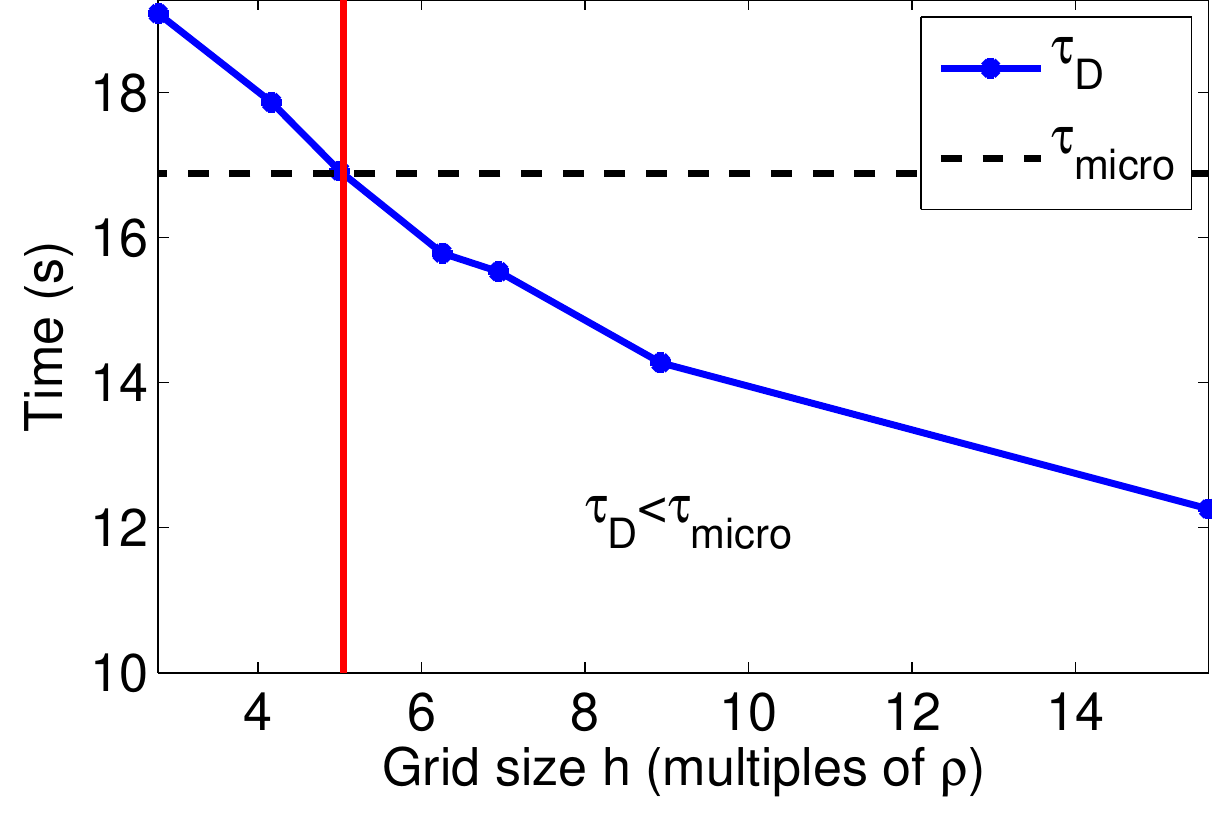}}
\subfigure[\, 3D]{\includegraphics[scale=0.5]{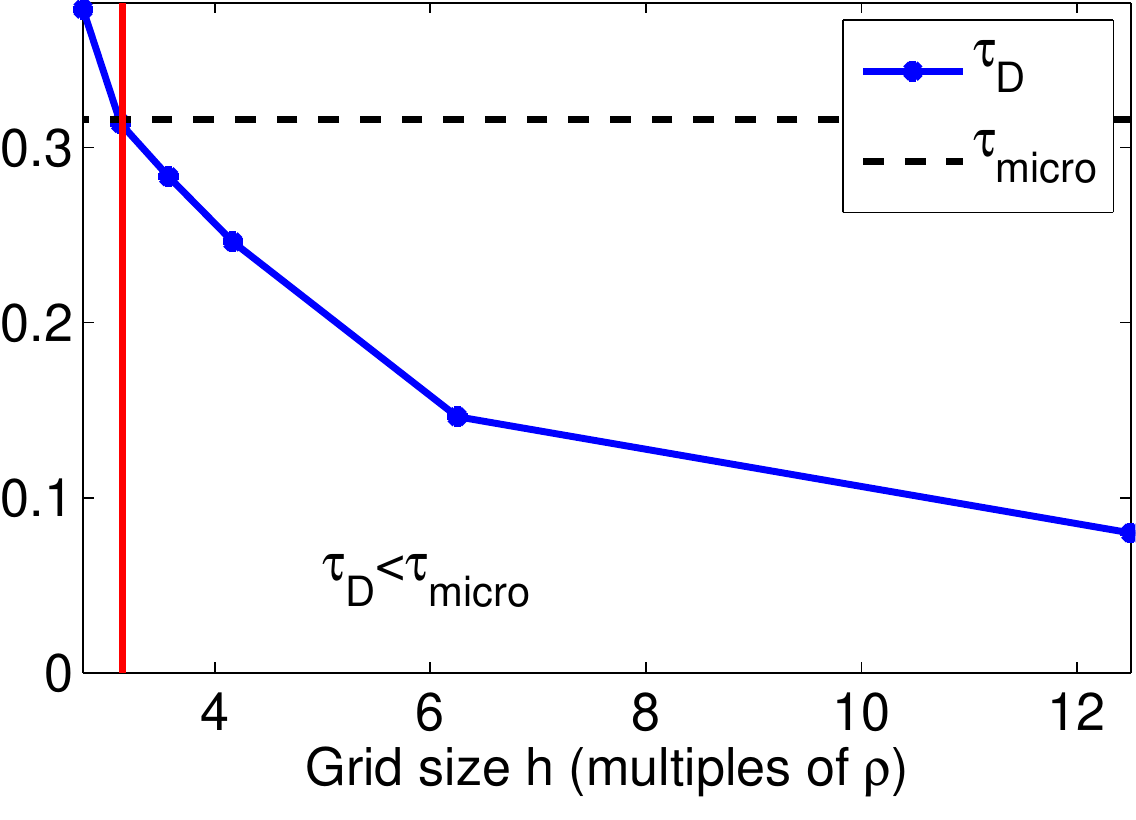}}
\caption{
The expected time until the molecules are in the same voxel for the first time, $\tau_D$, is computed with the RDME on a structured grid on a square (a) and a cube (b) with reflective boundary conditions. To the right of the vertical line we have $\tau_D<\tau_{\mathrm{micro}}$. In that region we could correct the mesoscopic reaction rate so that the expected time until the molecules react matches the expected time that we obtain by simulating the system on the microscale. To the left of the vertical line we have $\tau_D>\tau_{\mathrm{micro}}$, and it is impossible to change the mesoscopic reaction rate such that the result obtained with the RDME will match the result obtained from the Smoluchowski model. On the square we find that $h^*\approx 5.1\rho$ and on the cube $h^*\approx \pi\rho$. In our simulations we have used a uniform, structured Cartesian mesh and $\rho=2\cdot 10^{-9}m$, $D=10^{-12}m^2s^{-1}$ (3D)  and $D=10^{-14}m^2s^{-1}$ (2D).
}
\label{fig:fig1}
\end{figure}

Due to the potential computational advantages of using the RDME over the microscopic model, it is natural to try to correct the RDME to agree better with the Smoluchowski model for fine lattice spacings. The most obvious approach is to try to adjust the mesoscopic association rate constants in the RDME, making them dependent on the discretization. This would preserve the local nature of the reactions, and hence the low computational cost of the conventional RDME. One immediate consequence of the analysis in the previous section, however, is that below a certain mesh size $h^*$ no such local correction can make the mesoscopic mean association rate between two molecules agree with the microscopic model. In fact, for a given domain and model, this happens precisely  when $\tau_D>\tau_{micro}$. This is illustrated in Fig. \ref{fig:fig1} for the case of a square with side length $L=250\rho$ and a cubic domain with side length $L=100\rho$, with $k_r\to\infty$. 

If this critical mesh size always occurred for $h^*<\rho$, the problem would be less of a practical issue since the conventional RDME model makes little sense for mesh sizes smaller than the sum of the interaction radii of the molecules. However, in our example $h^*$ occurs for $h^*\approx \pi\rho$ (3D) and $h^*\approx 5.1\rho$ (2D). As long as $\tau_D<\tau_{micro}$, that is $h>h^*$, it is theoretically possible to modify the association rate, i.e. derive a discretization-dependent reaction propensity $q(h)$ that will give the same mean association time as the microscopic model. Clearly, to be able to give the correct binding times for as large a regime as possible, the modified reaction constant should have the property $q \to \infty$ for $\tau_D \to \tau_{micro}$.

While our analysis does not preclude the possibility to better match the mean association time by increasing $D$ and thus decreasing $\tau_D$, this would make the effective diffusion too fast and thus introduce another source of error. 

\section{Discussion}
\begin{figure}[htp]
\centering
\includegraphics[width=\linewidth]{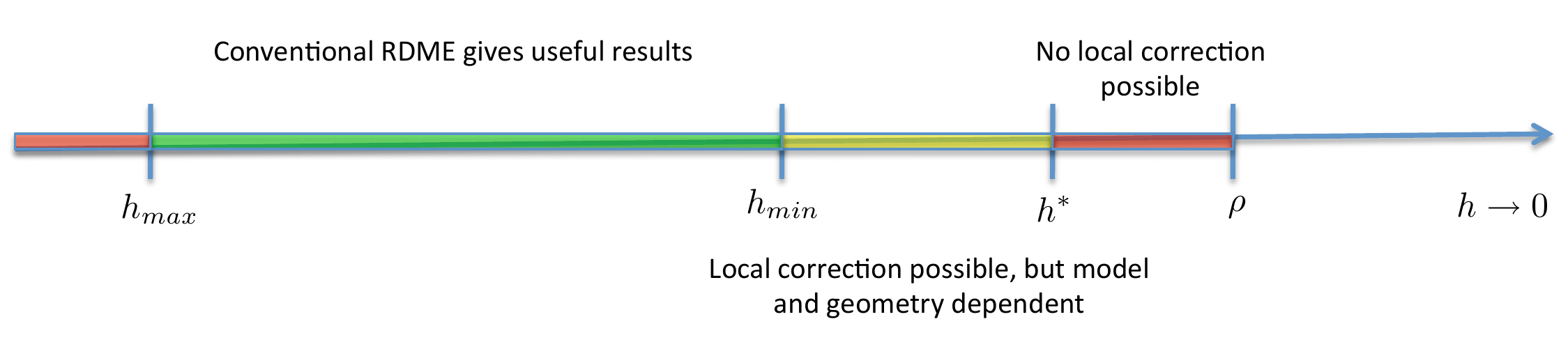}
\caption{
Schematic representation of the RDME's behavior as a function of $h$. For $h < h^*$, no local correction to the conventional mesoscopic reaction rates exists for the simple problem of diffusion to a target. 
}
\label{fig:hline}
\end{figure}

Fig. \ref{fig:hline} shows a schematic representation of the RDME's behavior as a function of the meshsize. For $h<\rho$, the RDME makes little sense and we cannot expect the model to work in this regime. For $h_{\mathrm{min}}<h<h_{\mathrm{max}}$ the conventional mesoscopic rate constants will work well, but for $h<h_{\mathrm{min}}$ the RDME will become increasingly inaccurate. For $h^* < h < h_{\mathrm{min}}$ it is possible to derive mesh and problem dependent reaction rates that make the RDME agree better with the microscopic model. 
The precise locations of the critical values $h_{\mathrm{min}}$, $h_{\mathrm{max}}$ and $h^*$ are model, geometry and discretization dependent.

In \cite{Isaacson2} Isaacson expands the analytical solutions of the RDME and the Smoluchowski equation for the simple model problem $A+B\to\emptyset$ in a series and computes the three leading terms in $h$. He shows that the two first terms will converge to the same value as $h$ tends to zero, but that the difference between the third term will diverge. There is an $h$ that minimizes the difference between the first terms of the expansion. This does not strictly prove that, for sufficiently small $h$, the error between the solutions increases as $h$ decreases, and the result is valid in an unbounded domain, but it still illustrates that for some $h<h_{\mathrm{min}}$ the reaction rates will need to be modified to make the mesoscopic model accurate. As we have shown here, however,  one will eventually reach $h^*$ and the difference between the models will inevitably increase.


Recently, two different corrections to $k_a$ have been proposed \cite{ErbanChapman,FBSE10}. In the 3D case with a cubic domain and a uniform, Cartesian discretization, Erban and Chapman \cite{ErbanChapman} consider the simple model problem $A+B\xrightarrow{k_a}B$, $\emptyset \xrightarrow{k_1} A$. They derive a mesh-dependent rate expression by matching the true steady-state distribution (which can be obtained analytically for this simple problem) to the distribution obtained using a meshsize $h$. They arrive at the following expression (in 3D)
\begin{equation}
q(h) = \frac{Dk_a}{Dh^3-\beta k_ah^2}.
\label{eq:propchap}
\end{equation}
They also find a critical mesh size $h_{crit}=\beta_{\infty}/({k_a}{D})$
under which no further correction can be made to $k_a$, where $\beta_{\infty} \approx 0.25272$ is a unitless constant valid for $L\gg h$. $h_{crit}$ makes the denominator in \eqref{eq:propchap} zero and hence $q \rightarrow \infty$ as $h \to h_{crit}$. Below that critical value of $h$, the mesoscopic association reaction is defined to $q=\infty$, i.e. association occurs as soon as the molecules are in the same voxel. Thus $q$ satisfies the basic requirements of our analysis: the existence of a critical meshsize and the correct limiting behavior as the meshsize tends to that critical value.  Substituting $k_a$ for the conventional expression and taking $k_r \to \infty$ we obtain $h_{crit} \approx \pi \rho$. Note that the propensity function \eqref{eq:propchap} was obtained without any comparison to the microscopic model. 

Fange et al. pursue a similar idea in \cite{FBSE10}. They study the problem $A+B\overset{k_a}{\underset{k_d}\leftrightharpoons} C$, and derive mesoscopic reaction rates such that the equilibration time of the system matches the equilibration time in the Smoluchowski model. They carry out this analysis in 2D and 3D on a disk and a sphere and obtain 
\begin{align*}
&p(h) = k_r/(1+\alpha(1-\beta)(1-0.58\beta)) \, \text{(3D)} \\
&p(h) = k_r/\{1+\alpha\ln[1+0.544(1-\beta)/\beta]\} \, \text{(2D)} 
\end{align*}
where $\beta = \rho/(\rho+\ell)$, $\rho+\ell$ is the radius of a disk with area $h^2$ in 2D and a sphere with volume $h^3$ in 3D, $\alpha = k_r/(4\pi\rho D)$ in 3D and $\alpha = k_r/(2\pi D)$ in 2D. 
These expressions do not predict a critical mesh size, but has the property $p(h) \to k_r~as~h\to 0$ in both 2D and 3D. 

Based on our analysis we easily obtain another correction in both 2D and 3D. Equation \eqref{eq:taumeso} suggest that we choose $k_{\mathrm{meso}}$ to be
\begin{align*}
k_{\mathrm{meso}} = (1+N_{\mathrm{steps}}^1)/(\tau_{\mathrm{micro}}-\tau_D)
\end{align*}
in order to have $\tau_{\mathrm{meso}} = \tau_{\mathrm{micro}}$. For $h$ sufficiently small we get
\begin{align}
\kme=\left\{\begin{array}{ll}
\frac{1+\left(L/h\right)^2}{\tau_{\mathrm{micro}}-[\frac{L^2}{2\pi D}\log\left(\frac{L}{h}\right)+\frac{0.1951L^2}{4D}]} & (2D)\\
\frac{1+(L/h)^3}{\tau_{\mathrm{micro}}-1.5164L^3/(6Dh)} & (3D)
\end{array}\right.
\label{eq:ourcorr}
\end{align}
From these expressions we obtain $h^* \approx \sqrt{\pi}\exp(0.1951\pi/2+3/4)\rho\approx 5.1\rho$ (2D) (where $\tau_{\mathrm{micro}}$ has been approximated by the analytical expressions for a disk derived in \cite{FBSE10}) and $h^* \approx \pi \rho$ (3D) (with $\tau_{\mathrm{micro}}$ approximated by $(k_a/L^3)^{-1}$) in good agreement with the simulations in Fig. 1.  

\begin{figure}[h!]
\centering
\subfigure[\, 2D]{\includegraphics[scale=0.5]{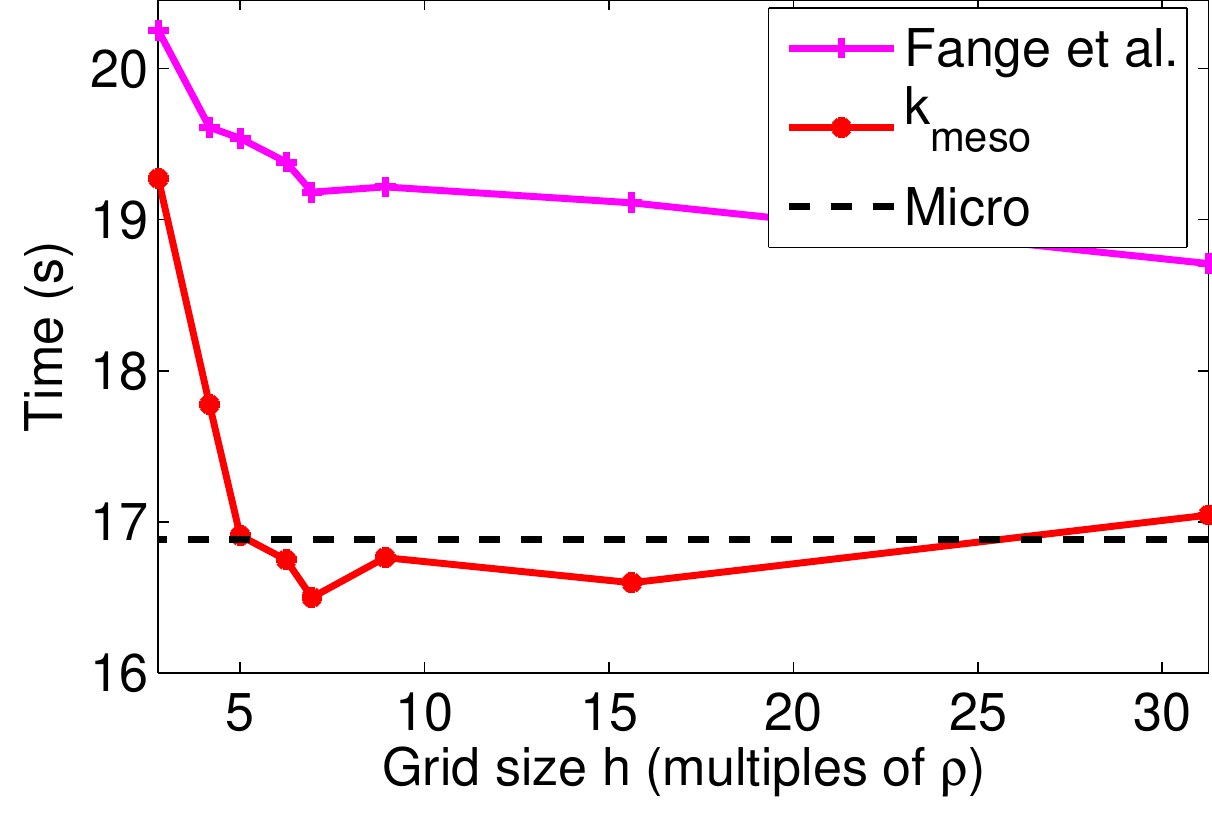}}
\subfigure[\, 3D]{\includegraphics[scale=0.5]{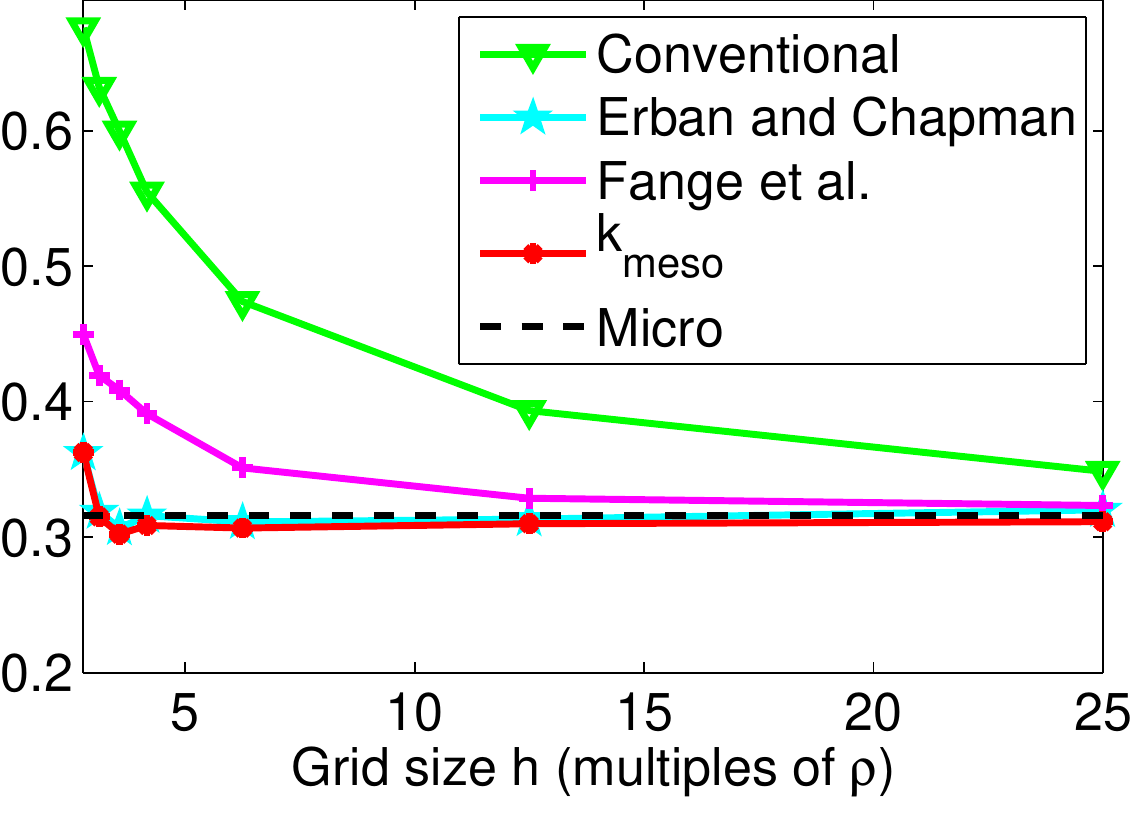}}
\caption{The expressions for the discretization-dependent mesocopic reaction rates from \cite{ErbanChapman,FBSE10} and those obtained here all depend on the ansatz used to derive them. For irreversible association in 2D (a), our correction give a different association time than that of Fange et al. and in 3D (b) our correction agrees well with that of Erban and Chapman but gives a different result than that of Fange et. al.  
Despite this, all expressions produce more accurate results than the conventional expression for our model problem for $h>h^*$.}
\label{fig:qa}
\end{figure}

The corrections obtained by Erban and Chapman do not coincide with the corrections obtained by Fange et al., illustrating how the corrections are dependent on the ansatz used to derive them. On the other hand, our corrections given by \eqref{eq:ourcorr} agree well with Erban and Chapmans in 3D and predict a similar $h^*$. 
As can be seen in Fig. \ref{fig:qa}, our corrections match the mean association time well in 2D (a) and all corrections give better results than the conventional expression $k_a$ in 3D (b). Interestingly, Fange et al. find experimentally for their example that they cannot match the Smoluchowski model in \cite[Fig. 3]{FBSE10} perfectly using the conventional RDME with the local corrections $p(h)$ even for $h\approx 5\rho$. Instead they modify the lattice model to allow for reactions between molecules in immediate neighboring voxels. In doing so they match the models all the way down to $h=2\rho$.  Our analysis explains why the local corrections alone were not sufficient, and their results demonstrate the possibility of better agreement with the microscopic model by a generalization of the conventional RDME to allow for neighbor-interactions.     

In this paper we have studied the behavior of the RDME for a simple but illustrative model problem: irreversible bimolecular association in the perfectly absorbing limit. For a finite but large $k_r$ (if the system is less diffusion-limited), $h^*$ would shift to the left towards $\rho$ in Fig. \ref{fig:fig1}. However, for a reversible reaction microscopic reversibility must hold and if $k_r$ is finite then the mesoscopic dissociation rate will be larger than the microscopic dissociation rate for some $h>h^*$. The $h^*$ that we obtain here in the irreversible case will thus be a lower bound for the reversible case. In conclusion, the conventional RDME cannot be made consistent with the Smoluschowski model since there will always be a meshsize for which no local correction to the reaction rate can give the correct mean association time. Above $h^*$ local corrections can extend the domain where the RDME works well. However, the corrections will inevitably be model and geometry dependent.

\section*{Acknowledgements} We acknowledge funding from the Swedish Research Council, the Royal Swedish Academy of Sciences FOA09H-63,
FOA09H-64, U.S DOE award DE-FG02-04ER25621 and Institute for Collaborative Biotechnologies Grant DAAD19-03-D-0004. 

\newcommand{\noopsort}[1]{}

\end{document}